\documentclass[tbtags,reqno]{amsart}

\usepackage{young}
  
\DeclareMathOperator{\Res}{Res}

\makeatletter
\renewcommand{\subsubsection}{\@startsection
{subsubsection}
{3}
{0mm}
{\baselineskip}
{-0.5\baselineskip}
{\normalfont\normalsize\bfseries}}
\makeatother

\newtheorem{theorem}{Theorem}
\newtheorem{lemma}[theorem]{Lemma}
\newtheorem{proposition}[theorem]{Proposition}

\newtheorem{formula}[theorem]{Formula}

\newtheorem{corollary}[theorem]{Corollary} 

\theoremstyle{definition}

\theoremstyle{remark}

\newtheorem*{acknow}{Acknowledgments}

\begin{document}

 \title[Rodrigues formulas for the Macdonald polynomials]{Rodrigues formulas for the Macdonald polynomials}
\author{Luc~Lapointe and Luc~Vinet}
\address
{Centre de Recherches Math\'ematiques, \\Universit\'e de
Montr\'eal,  C.P.~6128,  succursale~Centre-ville, \\ Montr\'eal,
Qu\'ebec, Canada, H3C 3J7}
 \date{}

\begin{abstract}
We present formulas of Rodrigues type giving the Macdonald polynomials for arbitrary partitions $\lambda$ through the repeated application of creation operators $B_k$, $k=1,\dots,\ell(\lambda)$ on the constant 1.  Three expressions for the creation operators are derived one from the other.  When the last of these expressions is used, the associated Rodrigues formula readily implies the integrality of the $(q,t)$-Kostka coefficients.  The proofs given in this paper rely on the connection between affine Hecke algebras and Macdonald polynomials.
\end{abstract}

\maketitle

\section{Introduction}
The Macdonald polynomials $J_{\lambda}(x;q,t)$ form a two-parameter basis for symmetric polynomials \cite{1}.  They play an important role in algebraic combinatorics and, in mathematical physics, they occur in particular in the wave functions of integrable quantum many-body models \cite{2}.  We shall show that the polynomials $J_{\lambda}(x)$ in $N$ variables can be constructed by acting with a string of creation operators $B_k$, $k=1,\dots,N$ on the constant 1, and shall thereby give Rodrigues formulas for these polynomials.  Such results were first obtained in \cite{3} in the limit case $q=t^{\alpha}$, $ t \to 1$ of the Jack polynomials and proved rather useful \cite{4}. 

 Three expressions $B_k^{(i)}$ $i=1,2,3$, will be obtained for the creation operators of the Macdonald polynomials.  Expression $B_k^{(1)}$ will first be derived using the Pieri formula.  The operator $B_k^{(2)}$ will then be shown to be equal to the operator $B_k^{(1)}$ and expression $B_k^{(3)}$ will finally be obtained from $B_k^{(2)}$ by observing that many terms in $B_k^{(2)}$ (and hence in $B_k^{(1)}$) act trivially on the Macdonald polynomials $J_{\lambda}$ associated to partitions with no more than $k$ parts.  Expression $B_k^{(1)}$ was first derived in \cite{5} where in addition, the $q$-difference operator version of  $B_k^{(3)}$ was given as a conjecture.  This third expression was also found by Kirillov and Noumi who provided two proofs \cite{6,7} of the fact that the operators $B_k^{(3)}$ are creation operators for the Macdonald polynomials.  It should be pointed out that the integrality of the $(q,t)$-Kostka coefficients \footnote{Other proofs of the integrality !
!
of the $(q,t)$-Kostka coefficien

ts have been given recently using different approaches by Garsia and Remmel \cite{11}, Garsia and Tesler \cite{12}, Knop \cite{13,14} and Sahi \cite{15}.} is an immediate consequence of the Rodrigues formula for $J_{\lambda}(x)$ associated to $B_k^{(3)}$.  We shall derive this formula from the one involving the operators $B_k^{(1)}$ by obtaining, as an intermediate step, the Rodrigues formula with the $B_k^{(2)}$ as creation operators.  Our proofs will rely in an essential way on the connection between affine Hecke algebras and Macdonald polynomials \cite{8,9}.  They will use in particular the fact that the Macdonald operators can be realized in terms of Dunkl-Cherednik operators.  This is the main difference between the proofs presented here and the short derivations given in \cite{10}.  It should be stressed that contrary to the approach followed in \cite{10}, the one taken here makes it possible to arrive at $B_k^{(2)}$ in a constructive fashion.

The outline of the paper is as follows.  In Section~2, basic facts about the relevant affine Hecke algebra realization are collected.  The Macdonald polynomials are introduced in Section~3 together with the commuting operators of which they are the simultaneous eigenfunctions.  Section~4 is the bulk of the paper.  This is where the expressions $B_k^{(1)}$, $B_k^{(2)}$ and $B_k^{(3)}$ are given (see (31),(32) and (34)) and derived (see Theorems 7,10 and 15).

\section{The affine Hecke algebra $H(\tilde W)$ \cite{6,8,9}}
  Let $\Lambda_N=\mathbb Q(q,t)[x_1,\dots,x_N]$ be the ring of polynomials in the $N$ variables $x_1,\dots,x_N$ with coefficients in $\mathbb Q(q,t)$, the field of rational functions in the two indeterminates $q$ and $t$.  The Weyl group $W \cong S_N$ is generated by the transpositions $s_i, i=1,\dots,N$.  On $x^{\lambda}=x_1^{\lambda_1}\dots x_N^{\lambda_N}$ their action is  such that
\begin{equation}
s_i x^{\lambda} = x^{s_i \lambda} s_i,
\end{equation}
where $s_i \lambda=(\lambda_{1},\dots,\lambda_{i-1},\lambda_{i+1},\lambda_{i},\lambda_{i+2}\dots ,\lambda_{N})$. We denote by $\Lambda_N^W$, the subring of all  symmetric polynomials. We can extend the action of the Weyl group $W$ on $\Lambda_N$ to one of  the affine Weyl group $\tilde W$ by introducing the elements $s_0$ and $\omega^{\pm 1}$ realized by:
\begin{equation}
\begin{split}
& s_0 = s_{N-1} \dots s_2 s_1 s_2 \dots s_{N-1} \tau_1 \tau_{N}^{-1}, \\
& \omega = s_{N-1} \dots s_1 \tau_1 = \tau_N s_{N-1} \dots s_1.
\end{split}
\end{equation}
where  $\tau_i$, the shift operator, is such that
\begin{equation}
\tau_i f(x_1,\dots,x_N) = f(x_1,\dots,q x_i,\dots,x_N)
\end{equation}
for any polynomial $f \in \Lambda_N $.

The generators of $\tilde W$ obey the fundamental relations:
\begin{equation}
\begin{aligned}
\mathrm{(i)} \ & s_i^2=1 \quad \quad \quad \quad \quad \quad &  & i=0,1,\dots,N-1,\\
\mathrm{(ii)} \ &  s_i s_j = s_j s_i  &      &   |i-j| \geq 2, \\
\mathrm{(iii)} \  & s_i s_j s_i=s_j s_i s_j &  &|i-j|=1, \\
\mathrm{(iv)} \ &  \omega s_i = s_{i-1} \omega  & &i=0,1,\dots,N-1.
\end{aligned}
\end{equation}
where the indices $0,1,\dots,N-1$ are understood as elements of $\mathbb Z_N = \mathbb Z/ N \mathbb Z$.
In the case of the  Weyl group $W$, for any $w \in W$, there is a smallest positive integer $p$ such that $w=s_{i_1}\dots s_{i_p}$ (reduced decomposition).  We say that $p$ is the length $L(w)$ of $w$.  Let $v,w \in W$, in the Bruhat order, $v \leq w$ if $v$ is of the form $v=s_{j_1}\dots s_{j_q}$ with $(j_1,\dots,j_q)$ a subsequence of $(i_1,\dots,i_p)$.

The operators
\begin{equation}
T_i = 1+ \frac{1-t^{-1}x_{i+1}/x_i}{1-x_{i+1}/x_i} (s_i -1),
\end{equation}
for $i=1,\dots,N-1$ and
\begin{equation}
T_0 = 1 + \frac{1-t^{-1}q^{-1}x_{1}/x_N}{1-q^{-1}x_{1}/x_N} (s_0 -1),
\end{equation}
and $\omega^{\pm 1}$ realize on $\Lambda_N$ the affine Hecke algebra $H(\tilde W)$ of $\tilde W$, that is they verify the defining relations
\begin{equation}
\begin{aligned}
\mathrm{(i)} \ & (T_i -1)(T_i+t^{-1})=0 \quad \quad \quad \quad & & i=0,1,\dots,N-1,\\
\mathrm{(ii)} \ & T_i T_j = T_j T_i & & |i-j| \geq 2,\\
\mathrm{(iii)} \ & T_i T_j T_i=T_j T_i T_j & & |i-j|=1, \\
\mathrm{(iv)} \ &  \omega T_i = T_{i-1}\omega  & & i=0,1,\dots,N-1,
\end{aligned}
\end{equation}
where again the indices are understood as elements of $\mathbb Z_N = \mathbb Z/ N \mathbb Z$.  The Dunkl-Cherednik operators $Y_1,\dots,Y_N$ are constructed as follows from the generators of $H(\tilde W)$:
\begin{equation}
Y_i = T_i \dots T_{N-1} \omega T_1^{-1} \dots T_{i-1}^{-1}.
\end{equation}
They form an Abelian algebra: $[Y_i,Y_j]=0$, $\forall i,j=1,\dots,N$. 
They also satisfy the following commutation relations with the $T_i$'s:
\begin{equation}
\begin{split}
& T_i Y_{i+1} T_i = Y_i, \\
& T_i Y_j = Y_j T_i  \quad \quad \quad j \neq i,i+1.
\end{split}
\end{equation}
 Let  $J = \{  j_1,  j_2, \dots , j_\ell \}$ denote  sets of cardinality $|J| = \ell$ made out  of integers $j_\kappa \in \{ 1, \dots , N \}$, $1 \le \kappa \le \ell$ such that $j_1 < j_2 < \dots < j_\ell$. We  introduce the operators
\begin{equation} 
Y_{J,u}
	= (1- u t^{\ell-1} Y_{j_1}) \dotsm (1- u t Y_{j_{\ell-1}})  (1- u Y_{j_{\ell}}),
\end{equation}
associated to such sets and labelled by a real number $u$.  If $|J|=0$, we define $Y_{J,u}=1$. 
 
For each $w=s_{i_1}\dots s_{i_p} \in W$, $T_w$ is defined by
\begin{equation}
T_w = T_{i_1}\dots T_{i_p}.
\end{equation}
Note that $T_w$ does not depend on the choice of the reduced decomposition of $w$.

The following relations between the generators of $H(\tilde W)$ and the variables $x_i$ will prove  useful
\begin{equation}
\begin{split}
& T_i x_{i} = x_{i+1} T_i - x_{i+1}(1-t^{-1}),\\
& T_i x_{i+1} = x_{i} T_i + x_{i+1}(1-t^{-1}),\\
& T_i x_{j} = x_{j} T_i  \quad \quad \quad j \neq i,i+1,\\
& T_i^{-1} x_{i} = x_{i+1} T_i^{-1} + x_{i}(1-t),\\
& T_i^{-1} x_{i+1} = x_{i} T_i^{-1} - x_{i}(1-t),\\
& \omega x_i = x_{i-1} \omega \quad \quad \quad i \neq 1, \\
& \omega x_1 = q x_N \omega.
\end{split}
\end{equation}
From (7) and (12), we see that the $x_i$'s, the $T_j$'s and $\omega^{\pm 1}$ form an algebra over the field $\mathbb Q(q,t)$.  An element $O$ of this algebra will be said to be normally ordered if all the variables $x_i$'s have been moved to the left, that is if $O$ has been put in the form
\begin{equation}
O=\sum_{\lambda}x^{\lambda}O_{\lambda},
\end{equation}
where $O_{\lambda}$ is in $\mathbb Q(q,t)[T_i,\omega^{\pm 1}]$.

\section{The Macdonald polynomials \cite{1}}
  Let $\lambda \in P \equiv \mathbb N^N$.  We denote by  $|\lambda|=\sum_{i} \lambda_i $, the degree of $\lambda$, and by $\ell(\lambda)$ the number of non-zero entries in $\lambda$.  The dominance order on the set $P^+ \subseteq P$ of all partitions $\lambda_1 \geq \lambda_2 \geq \dots \geq \lambda_N \geq 0$, is  $\lambda \ge \mu$ if $|\lambda|=|\mu|$ and $ \lambda_1 + \lambda_2 + \dots + \lambda_i \ge \mu_1 + \mu_2 + \dots + \mu_i$ for all $i$. This ordering is extended to $P$ as follows \cite{16}.  The orbit $W \lambda$ of $\lambda \in P$ under the action of the symmetric group $W \cong S_N$ will contain a unique partition $\lambda^+ \in P^+$.  We denote by $w_{\lambda}$, the unique element of minimal length such that $w_{\lambda} \lambda = \lambda^+$. We have $\lambda \ge \mu$ if either $\lambda^+ > \mu^+$ or $\lambda^+ = \mu^+$ and $w_{\lambda} \leq w_{\mu}$ in the Bruhat order of $W$.  Note that $\lambda^+$ is the unique maximum of $W \lambda$.

Homogeneous symmetric polynomials are labelled by partition $\lambda$ of their degree.  In the remainder of this section, $\lambda$  always stands for a partition, that is $\lambda \in P^+$.  

Three standard bases for $\Lambda_N^W$, the space of symmetric functions, are:
\begin{enumerate}
\item[(i)] the power sum symmetric functions $p_{\lambda}$ which in terms of the power sums 
\begin{equation}
p_i = \sum_k x_k^i,
\end{equation}
are given by
\begin{equation}
p_{\lambda} = p_{\lambda_1} p_{\lambda_2} \dots,
\end{equation}
\item[(ii)] the monomial symmetric functions $m_{\lambda}$ which are
\begin{equation}
m_{\lambda} = \sum_{\text{distinct permutations}} x_1^{\lambda_1} x_2^{\lambda_2} \dots
\end{equation}
\item[(iii)] the elementary symmetric functions $e_{\lambda}$  which in terms of the $i^{th}$ elementary function 
\begin{equation}
e_i = \sum_{j_1<j_2<\dots<j_i} x_{j_1} x_{j_2} \dots x_{j_i}= m_{(1^i)},
\end{equation}
are given by
\begin{equation}
e_{\lambda} = e_{\lambda_1} e_{\lambda_2} \dots .
\end{equation}
\end{enumerate}

The Macdonald polynomials can now be presented as follows.
To the partition $\lambda$ with $m_i(\lambda)$ parts equal to $i$, we associate the number
\begin{equation} \label{1}
z_\lambda
	= 1^{m_1} m_1 !  \, 2^{m_2} m_2! \dotsm
\end{equation}
We define a scalar product $\langle \ , \ \rangle_{q,t}$ on $\Lambda_N^W$  by
\begin{equation}
\langle p_\lambda, p_\mu \rangle_{q,t}
	=\delta_{\lambda \mu} z_\lambda \prod_{i=1}^{\ell(\lambda)} \frac{1-q^{\lambda_i}}{1-t^{\lambda_i}},
\end{equation}
where $\ell(\lambda)$ is the number of parts of $\lambda$.  The Macdonald polynomials $J_\lambda (x; q,t) \in \Lambda_N^W$ are uniquely specified by 
\begin{align}
\mathrm{(i)} \ &  \langle J_\lambda, J_\mu \rangle_{q,t} = 0, \qquad \text{if } \lambda \ne \mu, \label{3}\\
\mathrm{(ii)} \ &  J_\lambda = \sum_{\mu \le \lambda} v_{\lambda\mu}(q,t) m_\mu,  \label{4}\\
\mathrm{(iii)} \ &  v_{\lambda\lambda}(q,t)= c_{\lambda}(q,t),
\label{5}
\end{align}
where 
\begin{equation}
c_{\lambda}(q,t) = \prod_{s \in \lambda} (1-q^{a(s)} t^{\ell (s)+1}).
\end{equation}
As usual $a(s)$ and $\ell(s)$ denote the number of squares in the diagram associated to the partition $\lambda$  that are respectively to the south and east of the square $s$.

For $r = 1,\dots,N$, let $M_N^r$ denote the Macdonald operator
\begin{equation}
 \quad M_N^r=\sum_I t^{(N-r)r+r(r-1)/2} \tilde A_I (x;t) \prod_{i\in I} \tau_i,
\end{equation}
where the sum is over all $r$-element subsets $I$ of $\{1,\dots,N\}$,
\begin{equation}
\tilde A_I (x;t) =  \prod_{\begin{subarray}{c} i \in I \\ j  \in I^c 
\end{subarray}} \frac{ x_i - t^{-1} x_j}{x_i - x_j},
\end{equation}
and  $M_N^0 \equiv 1$. These operators commute with each other, $[M_N^r,M_N^l]=0$ and are diagonal on  the Macdonald polynomials basis.
From the Macdonald operators, one constructs
\begin{equation}
M_N(X;q,t) = \sum_{r=0}^N M_N^r X^r,
\end{equation}
with $X$ an arbitrary parameter. With $J$, a set of cardinality $|J|=j$, we shall also use $M_{J}(X;q,t)$ to represent the operator  $M_{j}(X;q,t)$ in the variables $x_i$, $i \in J$.  The generating function $M_N$ will play a crucial role in the following.  Its action on $J_{\lambda}(x;q,t)$ with $\ell(\lambda) \le N$ is given, remarkably, by
\begin{equation}
M_N(X;q,t) J_{\lambda}(x;q,t) = a_{\lambda}(X;q,t) J_{\lambda}(x;q,t),
\end{equation}
where 
\begin{equation}
a_{\lambda}(X;q,t) = \prod_{i=1}^N (1+X q^{\lambda_i} t^{N-i}).
\end{equation}
From (27) we see that the eigenvalue of $M_N^r$ on $J_{\lambda}(x;q,t)$ is the coefficient of $X^r$ in the polynomial (29).

It is known (see for instance \cite{6}) that the Macdonald operators can be rewritten in terms of Dunkl-Cherednik operators on $\Lambda_N^W$.  In particular, we have that
\begin{equation}
\Res Y_{\{1,\dots,N\},u}=M_N(-u;q,t),
\end{equation}
where $\Res O$ means that $O$ is restricted to act on $\Lambda_N^W$.

\section{Creation operators}
We now give the  expressions $B_k^{(i)}$ $i=1,2,3$; $k=1,\dots,N$ of the creation operators that we will derive in the remainder of the paper. 

{~~}

\noindent $\bullet$ Expression 1
\begin{equation}
B_k^{(1)} = \frac{1}{(q^{-1};t^{-1})_{N-k}} Y_{\{1,\dots,N\},t^{k+1-N}q^{-1}} e_k,
\end{equation}
where for $n$ positive integer, $(a;q)_n= (1-a)(1-qa) \dots (1-q^{n-1}a)$ and $(a;q)_0 \equiv 1$.
 $\bullet$ Expression 2
\begin{equation}
B_k^{(2)}=\sum_{|I|=k}x_I  \sum_{m=0}^{N-k} \sum_{\begin{subarray}{c} I' \subseteq I^c \\ |I'|=m \end{subarray}} \frac{q^{-m}}{(t^{k-N+1}q^{-1};t)_m} t^{-d(I',I^c)} Y_{I \cup I',t^{1-m}}.
\end{equation}
The quantity $d(I,J)$ entering in the above expression depends on nested subsets $J \subseteq I$ of $\{1,\dots,N\}$ and is defined as follows.  Order the elements of $I$ so that $I=\{i_0,\dots,i_{\ell-1}\}$ with $i_0<i_1< \dots <i_{\ell-1}$. Let $J =\{i_{j_1},\dots,i_{j_m}\} \subseteq I$ with its elements ordered so that $(i_{j_1},\dots,i_{j_m})$ is a subsequence of $(i_{0},\dots,i_{\ell-1})$, $0 \leq j_{\kappa} \leq \ell-1$; $k=1,\dots,m$. We then define 
\begin{equation}
d(J,I)= \sum_{k=1}^m j_k -m(m-1)/2.
\end{equation}
Note that the sum is over the indices that identify the elements of $J$ in the reference set $I$.  If $|J|=0$, then $d(J,I) \equiv 0$.

\noindent $\bullet$ Expression 3
\begin{equation}
B_k^{(3)} = \sum_{|I|=k} x_I Y_{I,t}.
\end{equation}

  We need to prove that these three sets of operators are such that 
\begin{equation}
B_k^{(i)} J_{\lambda}(x)=J_{\lambda+(1^k)}(x),
\end{equation}
if $\ell(\lambda) \leq k$.  If this is so, the following Rodrigues formula for the Macdonald polynomials associated to any partition $\lambda$ are easily seen to hold
\begin{equation}
J_\lambda(x;q,t)
	= (B_N^{(i)})^{\lambda_N} (B_{N-1}^{(i)})^{\lambda_{N-1} -\lambda_N} \dots (B_1^{(i)})^{\lambda_1 - \lambda_2} \cdot 1.
\end{equation}
That the  first expression has property (35) will follow from the Pieri formula which gives the action of the elementary symmetric functions $e_k$ on the monic Macdonald polynomials $P_{\lambda}=1/c_{\lambda} (q,t) J_{\lambda}$.  This formula reads \cite{1}
\begin{equation}
e_k P_{\lambda}= \sum_{\mu} \Psi_{\mu/\lambda} P_{\mu},
\end{equation}
where the sum is over all partitions $\mu$ containing $\lambda$ such that the set-theoretic difference $\mu-\lambda$ is $k$-dimensional with the property that $\mu_i -\lambda_i \leq 1$, $ \forall i \geq 1$.  With  $C_{\mu/\lambda}$  and  $R_{\mu/\lambda}$ respectively  denoting the union of the columns and of the rows  that intersect $\mu - \lambda$, the coefficients $\Psi_{\mu/\lambda}$ are given by
\begin{equation} 
\Psi_{\mu/\lambda}= \prod_{\begin{subarray}{c} s \in C_{\mu/\lambda} \\ s 
\not \in R_{\mu/\lambda} \end{subarray}} \frac{b_{\mu}(s)}{b_{\lambda}(s)}
\end{equation}
where
\begin{equation}
b_{\lambda}(s) = {\cases \frac{1-q^{a(s)} t^{\ell(s)+1}}{1-q^{a(s)+1} t^{\ell(s)}} \qquad & {\text{if}}~s \in \lambda \\
1    & {\text{otherwise}} \endcases}.
\end{equation}
We shall then construct $B_k^{(2)}$ from $B_k^{(1)}$, using the realization of the  affine Hecke algebra $H(\tilde W)$ given in Section~2 and properties of the Dunkl-Cherednik operators.  Upon proving the operator equality  $B_k^{(2)}=B_k^{(1)}$ we shall infer that $B_k^{(2)}$ are indeed creation operators.  Last, we shall prove to conclude the derivation that $B_k^{(3)}J_{\lambda}=$ $B_k^{(2)}J_{\lambda}=$ $J_{\lambda+(1^k)}$ on Macdonald polynomials with $\ell(\lambda) \leq k$. 
 
We start by giving  some results that we will need in the sequel.  First, a lemma that has to do with the normal ordering of some expressions (see(13)):
\begin{lemma}
\begin{equation}
   Y_k x_{\ell} = \sum_{j \geq \ell} x_j O_j,
\end{equation}
with $O_j \in \mathbb Q(q,t)[T_i,\omega^{\pm 1}]$ . And, 
\begin{equation}
Y_k x_1 \dots x_{\ell} =   
\cases
q x_1 \dots x_{\ell} Y_k \quad & {\text{if }} \ell \geq k\\
x_1 \dots x_{\ell} Y_k  & {\text{if }} \ell < k
\endcases
+ \sum_{|\lambda|=\ell} x^{\lambda} O_{\lambda},
\end{equation}
where all $\lambda$'s in the sum contain at least one non-zero part $\lambda_j$ with $j > \ell$, that is $\lambda \not \in P^+$.
\end{lemma}
It is easily proved by induction from (12).  A corollary of (40) is:
\begin{corollary}
For any $\lambda \in P, |\lambda|=\ell$, containing at least one non-zero part $\lambda_j$ with $j > \ell$, we have, for any $k$,
\begin{equation}
Y_k x^{\lambda}= \sum_{|\mu|=\ell} x^{\mu} O_{\mu},
\end{equation}
where all $\mu$'s in the sum contain  at least one non-zero part $\mu_j$ with  $j > \ell$, that is $\mu \not \in P^+$.
\end{corollary}
This is seen from (40) by commuting first $Y_k$  with one of the   $x_j$ with $j>\ell$ and $\lambda_j \neq 0$.

Next a lemma about the normal ordering of expressions involving $T_i$ and $x^{\lambda}$.
\begin{lemma} {~~~~~~~~~~~~~~~~~~~~~~~~~~~~~~~~~~~~~~~~~~~~~~~~~~~}

{~~~~~~~~~~~~}

\noindent (i) if $ s_i \lambda > \lambda$
\begin{equation}
T_i x^{\lambda}= x^{s_i \lambda} T_i + \sum_{\mu < s_i \lambda} x^{\mu} O_{\mu},
\end{equation}
(ii) if $ s_i \lambda < \lambda$
\begin{equation}
T_i x^{\lambda} = \sum_{\mu < \lambda}  x^{\mu} O_{\mu},
\end{equation}
(iii) if $s_i \lambda = \lambda$
\begin{equation}
T_i x^{\lambda}= x^{\lambda} T_i .
\end{equation}
\end{lemma}
Proof.  Since $T_i$ commutes with all the variables except $x_i$ and $x_{i+1}$, it suffices to look at the action of $T_i$ on $x_i^{\lambda_i}x_{i+1}^{\lambda_{i+1}}$.  The third case occurs when $\lambda_i=\lambda_{i+1}$ and it is trivially verified that $T_i (x_i x_{i+1})^{\lambda_i}= (x_i x_{i+1})^{\lambda_i} T_i$.  From this result, in case (i) where $\lambda_{i+1}>\lambda_i$, we see upon factoring $(x_i x_{i+1})^{\lambda_i}$ that it suffices to consider the action of $T_i$ on $x_{i+1}^{\lambda_{i+1}-\lambda_{i}}$.  Similarly, in case (ii), we see that we only need to consider how $T_i$ acts on   $x_i^{\lambda_i-\lambda_{i+1}}$. The proof is then  straightforwardly completed using (12).
\begin{lemma}
If $\mu$ and $\lambda$ with $\mu \neq \lambda$ are in the same orbit $W \lambda^+$ and such that $L(w_{\mu}) \geq L(w_{\lambda})\neq 0$ , then
\begin{equation}
T_{w_{\lambda}} x^{\mu}= \sum_{\rho < \lambda^+} x^{\rho} O_{\rho}.
\end{equation}
\end{lemma}
Proof.  The only non-trivial case is when $L(w_{\mu}) = L(w_{\lambda})$.  In this case, we have from Lemma~3 
\begin{equation}
T_{w_{\lambda}} x^{\mu}= T_{i_1} \dots T_{i_p} x^{\mu} = \sum_{\rho \leq w_{\mu,\lambda} \mu } x^{\rho} O_{\rho},
\end{equation}
with $w_{\mu,\lambda}$ some Weyl group element such that $w_{\mu,\lambda}< w_{\lambda}$ in the Bruhat order.  In order to have  $w_{\mu,\lambda}= w_{\lambda}$, case (i) of Lemma~3 would have to apply for every permutation $s_{i_k}$ in the reduction of $w_{\lambda}$, but this is impossible since it would require that  $ s_{i_{p-k}} ( s_{i_{p-k+1}} \dots s_{i_p} \mu) >  s_{i_{p-k+1}} \dots s_{i_p} \mu$ for   $k=1,\dots,p-1$, in other words, it would demand that  $w_{\lambda} \mu = \lambda^+$ which can not be  the case because $\mu \neq \lambda$ by hypothesis.  We thus have that all the $\rho$'s entering in (47) are such that $\rho \leq  w_{\mu,\lambda} \mu < \lambda^+$, which proves the lemma.
\begin{proposition}
If a non-zero operator is of the form $O= \sum_{|\mu|=k} x^{\mu} O_{\mu}$ with $O_{\mu}=0$ when $\mu \in P^+$,  there is at least one $T_i$, $i=1,\dots,N-1$, for which $T_i O \neq O T_i$ and hence $O$ is not symmetric.
\end{proposition}
Proof.  Suppose that $O$ is symmetric and of the form $O= \sum_{|\mu|=k} x^{\mu} O_{\mu}$ with $O_{\mu}=0$ when $\mu \in P^+$.  There exists one term $x^{\lambda} O_{\lambda}$ of $O$ with $O_{\lambda} \neq 0$ and  such that, either $\mu^+ \not > \lambda^+$ or  $ \mu^+ = \lambda^+$ and $L(w_{\mu}) \not < L(w_{\lambda})$, for all $\mu$ such that $O_{\mu} \neq 0$ in the decomposition of $O$.  From Lemma~3 (which imply that $T_w x^{\mu}=\sum_{\rho \leq \mu^+} x^{\rho}\tilde O_{\rho}$ for any $w \in W$ and $\mu \in P$) and Lemma~4, we then have that
\begin{equation}
T_{w_{\lambda}} O = x^{\lambda^+} T_{w_{\lambda}} O_{\lambda} + \sum_{\rho;\rho \neq \lambda^+} x^{\rho} O_{\rho}'.
\end{equation}
Since $O$ is assumed to be symmetric, we also have
\begin{equation}
T_{w_{\lambda}} O = O T_{w_{\lambda}} =\sum_{\mu} x^{\mu} O_{\mu}'',
\end{equation}
with $O_{\mu}''=0$ if $\mu \in P^+$.  From (48) and (49), since $T_{w_{\lambda}} O_{\lambda} \neq 0$, we have  a contradiction and hence Proposition~5 must be true.
\begin{corollary}
Two normally ordered symmetric operators $A$ and $B$, whose factors of  $x^{\lambda^+}$ are the same for any $  \lambda^+ \in P^+ $, are equal.
\end{corollary}
Since $A-B$ is symmetric and does not have any part in $x^{\lambda^+}$, $\forall \lambda^+ \in P^+$, it must be zero from Proposition~5.

\begin{theorem}
 For any partition $\lambda$ with $\ell(\lambda) \leq k$,
the operators $B^{(1)}_k$ act as follows on the Macdonald polynomials
$J_{\lambda}(x;q,t)$:
\begin{equation}
B^{(1)}_k J_{\lambda}(x) = J_{\lambda+(1^k)}(x).
\end{equation}
\end{theorem}
Proof.  The following lemma is an immediate consequence of the Pieri formula.
\begin{lemma} For $\lambda$ a partition such that $\ell(\lambda) \leq k$, the action of $e_k$ on $P_{\lambda}$ is given by 
\begin{equation}
e_k P_{\lambda} =  P_{\lambda+(1^k)} +
\sum_{\mu \neq \lambda+(1^k)} \Psi_{\mu/\lambda} P_{\mu},
\end{equation} 
where all the $\mu$'s in the sum are such that $\mu_{k+1}=1$.
\end{lemma}
Indeed, the only way to construct a $\mu$ with $\mu_{k+1} \neq 1$ is to add a 1 in each of the first $k$ entries of $\lambda$.
From Lemma~8 and (28) and (29), we have
\begin{equation}
 Y_{\{1,\dots,N\},t^{k+1-N}q^{-1}} e_k  P_{\lambda} = \prod_{i=1}^k (1-t^{k+1-i}q^{\lambda_i}) (q^{-1};t^{-1})_{N-k} P_{\lambda+(1^k)}
\end{equation}
since the eigenvalues
\begin{equation}
a_{\mu}(-t^{k+1-N}q^{-1};q,t) = \prod_{i=1}^N (1-t^{k+1-i}q^{\mu_i-1}),
\end{equation}
of  $Y_{\{1,\dots,N\},t^{k+1-N}q^{-1}}$ on the $P_{\mu}$'s in (51) vanish if $\mu_{k+1}=1$.

From the definition given in (24), it is easy to check that
\begin{equation}
\frac{c_{\lambda+(1^k)}}{c_{\lambda}}=\prod_{i=1}^k (1-t^{k+1-i}q^{\lambda_i}).
\end{equation}
Using this result and passing from $P_{\lambda}$ to $J_{\lambda}$ we see that 
\begin{equation}
B_k^{(1)}J_{\lambda}= \frac{1}{(q^{-1};t^{-1})_{N-k}} Y_{\{1,\dots,N\},t^{k+1-N}q^{-1}} e_k  J_{\lambda}= J_{\lambda+(1^k)}
\end{equation}
when $\ell(\lambda) \leq k$.  This proves Theorem~7.

We are going to order $B^{(1)}_k$ normally and thus  move all the variables contained in $e_k (x)$ to the left.  In doing so, we shall only focus on terms having  $x^{\lambda}$ with $\lambda \in P^+$ on the left, knowing from Corollary~6, that we only need to symmetrize these terms in order to find the full expression.  From Corollary~2, we see that the operators  $Y_k x_I$ will not have terms of the form $x^{\lambda^+}$ to the left whenever  $I \neq \{1,\dots,k \}$.  There thus only remains to consider
\begin{equation}
\begin{split}
\frac{1}{(q^{-1};t^{-1})_{N-k}} (1-t^k q^{-1} Y_1)& \dotsm (1- t q^{-1} Y_k)\\
\times &  (1- q^{-1} Y_{k+1}) \dotsm (1-t^{k+1-N} q^{-1} Y_N) x_1 \dots x_k .
\end{split}
\end{equation}
From (41) and Corollary~2, we see that the expression below is the only term of type  $x^{\lambda^+}$ in (56) once the variables $x_i$'s have been moved to the left:
\begin{equation}
\begin{split}
\frac{1}{(q^{-1};t^{-1})_{N-k}} x_1 \dots x_k (1-t^k Y_1) & \dotsm (1- t Y_k)\\
& \quad \times (1- q^{-1} Y_{k+1}) \dotsm (1-t^{k+1-N} q^{-1} Y_N).
\end{split}
\end{equation}
Thus,  
\begin{equation}
\frac{1}{(q^{-1};t^{-1})_{N-k}} x_1 \dots x_k Y_{\{1,\dots,k\},t} Y_{\{k+1,\dots,N\},t^{k+1-N}q^{-1}}
\end{equation}
is the only term of $B_k^{(1)}$ that has to the left a factor of the form $x^{\lambda^+}$.  Before symmetrizing, we shall expand this last expression using the following lemma:
\begin{lemma}
\begin{equation}
Y_{\{1,\dots,\ell\},t^{-\ell+1}q^{-1}}= \sum_{m=0}^{\ell} q^{-m}(q^{-1};t^{-1})_{\ell-m} \sum_{\begin{subarray}{c} I \subseteq \{1,\dots,\ell\} \\ |I|=m \end{subarray}} t^{-d \bigl( I,\{1,\dots,\ell \} \bigr) } Y_{I,t^{1-m}}.
\end{equation}
\end{lemma}
Proof.  We proceed by induction.  The formula is easily seen to hold in the case $\ell=1$.  Assuming that (59) is true, we thus have
\begin{equation}
\begin{split}
Y_{\{1,\dots,\ell+1\},t^{-\ell}q^{-1}}&=
Y_{\{1,\dots,\ell\},t^{-\ell+1}q^{-1}} (1-Y_{\ell+1}t^{-\ell}q^{-1}) \\
& = \sum_{m=0}^{\ell} q^{-m}(q^{-1};t^{-1})_{\ell-m} \Bigl(1-t^{-(\ell-m)}q^{-1}+t^{-(\ell-m)}q^{-1}\\
& \quad \quad \quad \quad \quad \quad \quad -Y_{\ell+1}t^{-\ell}q^{-1} \Bigr)  
 \sum_{\begin{subarray}{c} I \subseteq \{1,\dots,\ell\} \\ |I|=m \end{subarray}} t^{-d \bigl( I,\{1,\dots,\ell \} \bigr) } Y_{I,t^{1-m}}\\
& = \sum_{m=0}^{\ell} q^{-m}(q^{-1};t^{-1})_{\ell+1-m}  \sum_{\begin{subarray}{c} I \subseteq \{1,\dots,\ell+1\} \\ \ell+1 \not \in I; |I|=m \end{subarray}} t^{-d \bigl( I,\{1,\dots,\ell +1\} \bigr) } Y_{I,t^{1-m}}\\
& + \sum_{m=0}^{\ell} q^{-m-1}(q^{-1};t^{-1})_{\ell-m}  \sum_{\begin{subarray}{c} I \subseteq \{1,\dots,\ell+1\} \\ \ell+1 \in I; |I|=m+1 \end{subarray}} t^{-d \bigl( I,\{1,\dots,\ell +1\} \bigr) } Y_{I,t^{-m}}\\
& =\sum_{m=0}^{\ell+1} q^{-m}(q^{-1};t^{-1})_{\ell+1-m} \sum_{\begin{subarray}{c} I \subseteq \{1,\dots,\ell+1\} \\ |I|=m \end{subarray}} t^{-d \bigl( I,\{1,\dots,\ell+1 \} \bigr) } Y_{I,t^{1-m}},
\end{split}
\end{equation}
which concludes the proof.  In the derivation, we have  used the following two properties of the quantity $d(J,I)$:  $d \bigl( I,\{1,\dots,\ell\} \bigr) =d \bigl( I,\{1,\dots,\ell+1\} \bigr)$ if $\ell+1 \not \in I$ and  $d \bigl( I \cup \{\ell+1\},\{1,\dots,\ell+1\} \bigr) =d \bigl( I,\{1,\dots,\ell\} \bigr) +\ell-m$. 

With the help of Lemma~9 and of the identity $Y_{\{1,\dots,k\},t}Y_{I,t^{1-m}}=Y_{\{1,\dots,k\}\cup I,t^{1-m}}$ if $|I|=m$, expression (58) can be recast in the form 
\begin{equation}
x_1 \dots x_k \sum_{m=0}^{N-k} \sum_{\begin{subarray}{c} I \subseteq \{k+1,\dots,N\} \\ |I|=m \end{subarray}} \frac{q^{-m}}{(t^{k-N+1}q^{-1};t)_m} t^{-d \bigl( I,\{k+1,\dots,N\} \bigr) } Y_{\{1,\dots,k\} \cup I,t^{1-m}}.
\end{equation}

We shall now give an expression in normal order (see (62)) which has (61) as its only term of type $x^{\lambda^+}$.  By Corollary~6, in order to show that this expression coincides with $B_k^{(1)}$ we shall only need to prove that it is symmetric.
\begin{theorem}
\begin{equation}
B_k^{(1)}=B_k^{(2)} \equiv \sum_{|I|=k}x_I  \sum_{m=0}^{N-k} \sum_{\begin{subarray}{c} I' \subseteq I^c \\ |I'|=m \end{subarray}} \frac{q^{-m}}{(t^{k-N+1}q^{-1};t)_m} t^{-d(I',I^c)} Y_{I \cup I',t^{1-m}},
\end{equation}
hence $B_k^{(2)}$ is also such that $B_k^{(2)}J_{\lambda}=J_{\lambda+(1^k)}$ for $\ell(\lambda) \leq k$.
\end{theorem}
Proof.  We first give some expressions that are easily checked  to commute with $T_i$ from properties (7),(9) and (12).  With $f \in \Lambda_N^W$,
\begin{equation}
\begin{split}
\mathrm{(I)} \ & (T_i-1) (x_i+x_{i+1})f=0 \\
\mathrm{(II)} \ & (T_i-1) x_i x_{i+1}f=0 \\
\mathrm{(III)} \  & (T_i-1) (1-u t Y_i) (1-u Y_{i+1})f=0 \\
\mathrm{(IV)} \ & (T_i-1) \Bigl[ x_i (1-u  Y_i)+x_{i+1} (1-u Y_{i+1}) \Bigr]f=0 \\
\mathrm{(V)} \ & (T_i-1) \Bigl[  (1-u  Y_i)+ t^{-1} (1-u Y_{i+1}) \Bigr]f=0 
\end{split}
\end{equation}
That $B_k^{(2)}$ has expression (61) has its only term of type $x^{\lambda^+}$ is obvious.  In view of the remark made before Theorem~10, we now only need to verify that the operators  $B_k^{(1)}$ as defined in (32) (and (62)) are symmetric, in other words that they satisfy $(T_i -1 )B_k^{(2)}f=0$, $\forall i=1,\dots,N-1$.  Since two sets, $I$ and $I'$, enter in the definition of $B_k^{(2)}$ we proceed by looking at all the inclusion possibilities of the indices $i$ and $i+1$ into these two sets and then examine the particular terms in $B_k^{(2)}$ that are affected by the action of $T_i$.  What we find using (63) is that these terms are separately or pairwise symmetric.  Indeed consider the cases: 
\begin{enumerate}
\item[(i)] $i,i+1 \not \in I$ and $i,i+1 \not \in I'$: $T_i$  trivially commutes.
\item[(ii)] $i,i+1 \not \in I$ and $i \in I',i+1 \not \in I'$: add the case $\bar I=I,\bar I'=(I' \cup \{i+1\}) \backslash \{i\}$ which is such that $ d(\bar I',\bar I^c)=d( I', I^c)+1$.  The symmetry follows from $\mathrm{(V)}$.
\item[(iii)] $i,i+1 \not \in I$ and $i,i+1  \in I'$: $T_i$ commutes owing to $\mathrm{(III)}$.
\item[(iv)] $i \in I,i+1 \not \in I $ and $i+1 \not  \in I'$:  add the case  $\bar I=(I \cup \{i+1\}) \backslash \{i\},\bar I'=I'$ which is such that $ d(\bar I',\bar I^c)=d( I', I^c)$. The symmetry is verified with the help of $\mathrm{(IV)}$.
\item[(v)] $i \in I,i+1 \not \in I$ and $i+1 \in I'$:  add the case $\bar I=(I \cup \{i+1\}) \backslash \{i\},\bar I'=(I' \cup \{i\}) \backslash \{i+1\}$ which is such that $ d(\bar I',\bar I^c)=d( I', I^c)$.  The symmetry is then confirmed using $\mathrm{(I)}$ and $\mathrm{(III)}$ .
\end{enumerate}
All the other cases are immediate consequences of  these cases and of (II).  Theorem~10 is thus seen to hold.

When going from $B_k^{(2)}$ to $B_k^{(3)}$, it is useful to obtain $\Res B_k^{(2)}$, that is the $q$-difference operator version of $B_k^{(2)}$. The next lemma gives the  essential step.
\begin{lemma}
\begin{equation}
\Res \sum_{|I|=k}x_I   \sum_{\begin{subarray}{c} I' \subseteq I^c \\ |I'|=m \end{subarray}}  t^{-d(I',I^c)} Y_{I \cup I',t^{1-m}}
= \sum_{|I|=k}x_I   \sum_{\begin{subarray}{c} I' \subseteq I^c \\ |I'|=m \end{subarray}}  \tilde A_{I \cup I'} M_{I \cup I'}(-t^{1-m};q,t).
\end{equation}
\end{lemma}
Proof.  We first give two formulas that we will need.  The first one is a well known identity and the second is a special case of a formula given by Garsia and Tesler (\cite{12}, Proposition~3.1).
\begin{formula}
\begin{equation}
e_m(1,t^{-1},\dots,t^{-N+1})= t^{m(m-1)/2} t^{-m(N-1)}
\begin{bmatrix}
N \\
m
\end{bmatrix}_t ,
\end{equation}
where the $q$-binomial coefficient is
\begin{equation}
\begin{bmatrix}
n\\
k
\end{bmatrix}_q
=\frac{(q;q)_{n}}{(q;q)_{k}(q;q)_{n-k}}.
\end{equation}
\end{formula}
\begin{formula}
With $N=|J \cup J^c|$, we have that, for any $k=0,\dots,N$ and $m=0,\dots,N-k$,
\begin{equation}
\sum_{|J|=k} x_J \sum_{\begin{subarray}{c} J' \subseteq J^c \\ |J'|=m \end{subarray}}  \tilde A_{J \cup J'}=t^{-m(N-k-m)} 
\begin{bmatrix}
N-k \\
m
\end{bmatrix}_t
 \sum_{|J|=k} x_J.
\end{equation}
\end{formula}
Given that 
\begin{equation}
 \sum_{\begin{subarray}{c}  |J'|=m \\ J' \subseteq J^c  \end{subarray}} t^{-d(J',J^c)} =
t^{m(m-1)/2} e_m(1,t^{-1},\dots,t^{-(N-k)+1}),
\end{equation}
for any subset $J$ of $\{1,\dots,N\}$ of cardinality $k$, 
the following lemma follows from Formulas~12 and 13.
\begin{lemma}
With $N=|J \cup J^c|$, we have that, for any $k=0,\dots,N$ and $m=0,\dots,N-k$,
\begin{equation}
 \sum_{|J|= k}  x_J \sum_{\begin{subarray}{c}  |J'|=m \\ J' \subseteq J^c  \end{subarray}} t^{-d(J',J^c)} =
\sum_{|J|= k}  x_J \sum_{\begin{subarray}{c}  |J'|=m \\ J' \subseteq J^c  \end{subarray}} \tilde A_{J \cup J'}.
\end{equation}
\end{lemma}
We now return to the proof of (64).  Since both sides  are symmetric, it will suffice to show that the coefficients of $\tau_1 \dots \tau_{\ell}$ for $\ell \leq N$ are identical on both sides of the equation.

To that end, let 
\begin{equation}
\begin{split}
& I=L \cup J, \quad \quad I'=\bar L \cup J',\\
& L \subseteq \{1,\dots, \ell\}, \quad \quad \bar L = \{1,\dots,\ell\} \backslash L,\\
& J,J' \subseteq \{1,\dots,N\} \backslash \{1, \dots,\ell \}= J \cup \bar J,\\
& J \cap \bar J= \phi, J' \subseteq \bar J,
\end{split}
\end{equation}
and define
\begin{equation}
[\ell,k]= {\cases \ell \qquad & {\text{if}}~\ell \leq k \\
k    &  {\text{if}}~\ell > k \endcases}.
\end{equation}

The only place in the l.h.s. of (64) where the operator product $\tau_1 \dots \tau_{\ell}$ will occur is in $\Res Y_1 \dots Y_{\ell}$ (see \cite{6} Proposition 6.1 and formula 7.20).  The coefficient of $\tau_1 \dots \tau_{\ell}$ in this expression is  $\tilde A_{\{1,\dots,\ell \}}$.  With this knowledge and the help of  (10),(25) and (27), we find that the coefficients of $\tau_1 \dots \tau_{\ell}$ on both sides of (64) are respectively:
\begin{equation}
\begin{split}
{\text{l.h.s.}} |_{\tau_1 \dots \tau_{\ell}}= t^{\ell(\ell-1)/2} \bigl( -t^{k-\ell+1} & \bigr)^{\ell} \tilde A_{\{1,\dots,\ell\}} \\
& \times \sum_{n=0}^{[\ell,k]}   \sum_{|L|=n}  x_L \sum_{|J|=k-n} x_J   \sum_{|\bar L \cup J'|=m} t^{-d(J',\bar J)}
\end{split}
\end{equation}
and 
\begin{equation}
\begin{split}
{\text{r.h.s.}} |_{\tau_1 \dots \tau_{\ell}}= t^{\ell(\ell-1)/2} \bigl( -t^{1-m}\bigr)^{\ell} & t^{\ell(m+  k-\ell)} \tilde A_{\{1,\dots,\ell\}} \\ 
& \times \sum_{n=0}^{[\ell,k]}   \sum_{|L|=n}  x_L \sum_{|J|=k-n} x_J   \sum_{|\bar L \cup J'|=m} \tilde A_{J \cup J'}^{J \cup \bar J},
\end{split}
\end{equation}
with
\begin{equation}
\tilde A_J^I = \prod_{\begin{subarray}{c} i \in J \\ j \in I \backslash J \end{subarray}} \frac{x_i - t^{-1}x_j}{x_i - x_j}.
\end{equation}
We have used in (72) the fact that $d(I',I^c)=d(J',\bar J)$ and in (73), the identity $\tilde A_{I \cup I'} \tilde A_{\{1,\dots,\ell\}}^{I \cup I'}= \tilde A_{\{1,\dots,\ell \}} \tilde A_{I \cup I' \backslash \{1,\dots,\ell\}}^{\{1,\dots,N\} \backslash \{1,\dots,\ell\}}$.  It is then immediate to see that the equality
\begin{equation}
{\text{l.h.s. of (64)}}|_{\tau_1 \dots \tau_{\ell}}={\text{r.h.s. of (64)}} |_{\tau_1 \dots \tau_{\ell}}
\end{equation}
holds, since after trivial simplifications it is seen to amount to
\begin{equation}
\begin{split}
& \sum_{n=0}^{[\ell,k]}   \sum_{|L|=n}  x_L \Bigl( \sum_{|J|=k-n} x_J   \sum_{\begin{subarray}{c}  |J'|=m-\ell+n \\ J' \subseteq \bar J  \end{subarray}}  t^{-d(J',\bar J)} \Bigr)=\\
& \quad \quad \quad \quad \quad \quad \quad \quad \quad \sum_{n=0}^{[\ell,k]}   \sum_{|L|=n}  x_L \Bigl( \sum_{|J|=k-n} x_J  \sum_{\begin{subarray}{c}  |J'|=m-\ell+n \\ J' \subseteq \bar J  \end{subarray}}  \tilde A_{J \cup J'}^{J \cup \bar J} \Bigr)
\end{split}
\end{equation}
and hence to follow from Lemma~14.

Once Lemma~11 is proved, the connection between $B_k^{(2)}$ and $B_k^{(3)}$ is readily obtained.
\begin{theorem}
For any partition $\lambda$, such that $\ell (\lambda) \leq k$, the actions of $B_k^{(2)}$ and $B_k^{(3)}$ on the Macdonald polynomials $J_{\lambda}(x)$ coincide:
\begin{equation}
B_k^{(3)}J_{\lambda}(x)=  B_k^{(2)} J_{\lambda}(x)= J_{\lambda+(1^k)}(x).
\end{equation}
\end{theorem}
This is shown to be true with the help of the following lemma
\begin{lemma} 
Let  $|I|=k$ and $|I'|=m$, $I' \subseteq I^c$.
\begin{equation}
M_{I \cup I'} (-t^{1-m};q,t) J_{\lambda}(x;q,t)=0,
\end{equation}
if $\ell(\lambda) \leq k $ and  $m > 0$.
\end{lemma}
Proof.  Denote by $x(I)$  the set of variables $\{ x_i, i \in J\}$.  The Macdonald polynomials are known \cite{1} to enjoy the property according to which
\begin{equation}
J_{\lambda} \bigl( x(I),x(I^c) \bigr)  = \sum_{\mu,\nu} \tilde f_{\mu \nu}^{\lambda} J_{\mu} \bigl( x(I) \bigr) J_{\nu} \bigl( x(I^c) \bigr)
\end{equation}
with $ \tilde f_{\mu \nu}^{\lambda}=0$ unless $ \mu \subset \lambda$ and  $ \nu \subset \lambda$ and in particular if $\ell(\mu)$ or $\ell(\nu)$ is greater than $k$. 

  Since $M_{I \cup I'} (-t^{1-m};q,t)$ is a $q$-difference operator depending only of  the variables $x_i$, $i \in I \cup I'$, we see from (79) that 
\begin{equation}
M_{I \cup I'} J_{\lambda}(x)  = \sum_{\mu,\nu} \tilde f_{\mu \nu}^{\lambda} J_{\mu} \bigl( x((I \cup I')^c) \bigr) M_{I \cup I'} J_{\nu} \bigl( x(I \cup I') \bigr).
\end{equation}
The proof of Lemma~16 is then completed by observing from (28)-(29) that
\begin{equation}
M_{I \cup I'} (-t^{1-m};q,t)  J_{\nu} \bigl( x(I \cup I') \bigr) = \prod_{i=1}^{k+m} (1-q^{\nu_i}t^{k+1-i}) J_{\nu} \bigl( x(I \cup I') \bigr)=0,
\end{equation}
whenever  $m>0$, since $\nu_{k+1}=0$.

Theorem~15 is thus an immediate consequence of this lemma since, by (64), we have that:
\begin{equation}
\Bigl( \sum_{|I|=k}x_I   \sum_{\begin{subarray}{c} I' \subseteq I^c \\ |I'|=m \end{subarray}}  t^{-d(I',I^c)} Y_{I \cup I',t^{1-m}} \Bigr) J_{\lambda}(x;q,t)=0,
\end{equation}
if $\ell(\lambda) \leq k$ and $m>0$;  this only leaves the $m=0$ term of $B_k^{(2)}$, which coincides with $B_k^{(3)}$.

\section{Conclusion}

 As already mentioned in the introduction,  the integrality of the $(q,t)$-Kostka coefficients follows quite straightforwardly from the Rodrigues formula (36) when the operators $B_k^{(3)}$ are used as creation operators.  This is explained  in \cite{6,7,10}.  Other interesting properties of $B_k^{(3)}$ have been conjectured in \cite{5}.  In particular, a formula that would give the action of these operators on arbitrary Macdonald polynomials has been proposed:  it looks like a deformation of the Pieri formula and, if true, would imply that the $N$ operators $ \Res \sum_{|I|=m} x_I Y_{I,t^{\kappa-m+1}}$, $m=1,\dots,N$  form a commuting set  for any $\kappa \in \mathbb R$.  We hope that the constructive approach presented in this paper will allow one to make progress towards proving these conjectures and unravelling the  algebra of which the creation operators are part of.

\begin{acknow}
We would like to express our thanks to Adriano Garsia  for various comments and suggestions.  

\noindent This work has been supported in part through funds provided by NSERC (Canada) and FCAR (Qu\'ebec).  L.~Lapointe holds a NSERC postgraduate scholarship.
\end{acknow}

\end{document}